# Notes on Feige's gumball machines problem

John H. Elton
Georgia Institute of Technology

**Abstract**
Uriel Feige posed this problem, as communicated by Peter Winkler in the puzzle column of the Communications of the ACM in August 2009: For $n$ non-negative integer valued independent random variables, each having mean one, what is the maximum probability that the sum is greater than $n$? In this note we reduce the problem to the case of two-valued random variables, and give a rigorous and detailed proof of the conjectured solution in the special case that the random variables are assumed identically distributed.

## Introduction

Let $X_i, i = 1,...n$ be independent random variables with ranges in the set of non-negative integers, and each with expected value 1, but not necessarily identically distributed. Let $L$ be the infimum of $P(\sum_{i=1}^{n} X_i \leq n)$ over all possible distributions satisfying the above conditions, where $P(E)$ is the probability of event E. $L$ is attained by some distribution and so is a minimum. The goal is to find $L$, and the distribution which attains it.

This problem is due to Uriel Feige, and was exposed as a puzzle in Peter Winkler's puzzle column in the Communications of the ACM, August 2009 [2, pg 104-105]. The random variables there were described as the output of $n$ independent gumball machines, with the subject playing each machine exactly once and counting the total number of gumballs obtained. The problem was stated as asking what distribution gives the maximum probability of getting strictly more than $n$ gumballs. The conjecture is that the optimal solution is to put probability $1/(n+1)$ for getting $n+1$ gumballs, and the remaining probability on zero gumballs.

In this note, we show how to reduce the problem to one of distributions on only two points, and we give a rigorous and detailed solution for the case that the machines are known to have the same distribution, which already requires quite a careful analysis of the numerics of the binomial distribution. This gives some insight into the difficulties that can be expected in the general case where the random variables are not assumed to be identically distributed.

Our reduction to two point distributions makes it possible to verify experimentally that the solution to the general case is the conjectured one up to $n = 20$ (we could go further but there seems no reason to at this time).

## Reduction to two-valued random variables

Let $p_i(j) = P(X_i = j), i = 1,...,n$, be the probability mass functions of the random variables. Assume throughout that $n \geq 2$.

**Lemma 1**. For a distribution obtaining the minimum, $p_i(j) = 0$ for $j > n+1, i = 0,...,n$.

Proof: Suppose the minimum is attained for some distribution with mass functions $p_i, i = 1,...,n$. Note that the expected value constraint along with the properties of a probability mass function implies $p_i(0) = \sum_{j \geq 2}(j-1)p_i(j)$ and $p_i(1) = 1 - \sum_{j \geq 2} jp_i(j)$.

Suppose $p_m(k) > 0$ for some $m$ and some $k > n+1$. Define

$$p'_m(k) = 0, \; p'_m(n+1) = p_m(n+1) + \frac{k}{n+1}p_m(k), \; p'_m(0) = p_m(0) - \left(\frac{k}{n+1} - 1\right)p_m(k),$$

and $p'_m(j) = p_m(j)$ for all other $j$. Note that $p'_m(0) > 0$ because $p_m(0) \geq (k-1)p_m(k)$.

Replacing the mass function $p_m$ for $X_m$ by $p'_m$, we see that $P(\sum_{i=1}^n X_i \leq n)$ is decreased because $P(X_m = 0)$ has decreased, $P(X_m = j)$ is unchanged for $1 \leq j \leq n$, and

$$P(\sum_{i=1}^n X_i \leq n) = P(X_m = 0)P(\sum_{i \neq m} X_i \leq n) + \sum_{j=1}^n P(X_m = j)P(\sum_{i \neq m} X_i \leq n - j).$$ 

This contradicts the hypothesis of minimality, so there can be no mass on values greater than $n+1$. □

**Lemma 2**. For a distribution obtaining the minimum, there exists $2 \leq j_i \leq n+1, i = 1,...,n$ such that $p_i(j) = 0$ for $j \neq 0$ and $j \neq j_i, i = 1,...,n$; thus $p_i(0) = 1 - 1/j_i, p_i(j_i) = 1/j_i$, and $p_i(j) = 0$ for all other $j$. So each random random variable has its mass concentrated on just two integers, one of which is zero and the other greater than or equal 2 and less than or equal $n+1$.

Proof: Write $P(\sum_{i=1}^n X_i \leq n) = \sum_{j=1}^n P(\sum_{i=0}^{n-1} X_i \leq n - j)P(X_n = j) = \sum_{j=0}^n c_j y_j$ where

$c_j = P(\sum_{i=0}^{n-1} X_i \leq n - j)$ and $y_j = p_n(j), j = 0,...,n+1$. Thus $P(\sum_{i=1}^n X_i \leq n)$ is a linear function of the $y_j$ considering the distribution of $X_1,..., X_{n-1}$ as fixed. As mentioned above, from the expected value constraint and the probability mass function constraint, we may obtain $y_0 = \sum_{j=2}^{n+1}(j-1)y_j, y_1 = 1 - \sum_{j=2}^{n+1} jy_j$, so we may write

$P(\sum_{i=1}^n X_i \leq n) = f(y_2, y_3,...y_{n+1})$ as a linear function of the variables $y_2,..., y_n$ subject to the linear inequalities $y_j \geq 0, j = 2,...,n+1$ and $\sum_{j=2}^{n+1} jy_j \leq 1$. So the domain S of the function f is a convex set, in fact the simplex bounded by the coordinate hyperplanes together with the hyperplane $\sum_{j=2}^{n+1} jy_j \leq 1$. Since f is linear, its minimum occurs at an extreme point of S, so either $y_j = 0, j = 2,...,n+1$ and so $y_1 = 1$ and $y_0 = 0$, or else $y_k \neq 0$ for some $2 \leq k \leq n+1$, $y_j = 0$ for $j \neq k, 2 \leq k \leq n+1$, and $ky_k = 1$. In this latter case, $y_1 = 0, y_k = 1/k,$ and $y_0 = 1 - 1/k$.

The same argument may be made for $X_i, i = 1,...,n-1$, which means that for each random variable the mass function is concentrated on at most two points, for a distribution attaining the minimum. We can eliminate the single point support on 1 and will do that later sometime. □

Remark: This lemma makes it easy to check the conjecture for small $n$ on a computer, and we have done that up to $n = 20$.

**The case that the random variables are known to be identically distributed.**

For now, let's get the result we seek with the additional assumption that the random variables are identically distributed. In that case, there exists $2 \leq j \leq n+1$ such that $p_i(0) = 1 - \frac{1}{j}, p_i(j) = \frac{1}{j}$. Then $P(\sum_{i=1}^{n} X_i \leq n) = \sum_{kj \leq n} \binom{n}{k} \left(\frac{1}{j}\right)^k \left(1 - \frac{1}{j}\right)^{n-k}$. We want to show this is minimum when $j = n+1$. We look at the behavior as a function of $j$ for $2 \leq j \leq n$, and it is a little stranger than you might expect. We found that it was easier to look first at cumulative binomial distributions with integer means, with the probability of success in a trial being $m/n$, and to use that to put bounds on the case when the probabilities are $1/j$.

**Lemma 3**. $f(n,m) := \sum_{k \leq m-1} \binom{n}{k} \left(\frac{m}{n}\right)^k \left(1 - \frac{m}{n}\right)^{n-k} \geq \frac{3}{8}$ for $n \geq 100, 12 \leq m \leq \frac{n}{2} + 1$

Proof: The sum is the probability that a binomial random variable with mean $m$ is less than or equal $m-1$. By a result of Kaas and Buhrman [1], the median is also $m$, so this exceeds $\frac{1}{2} - \binom{n}{m} \left(\frac{m}{n}\right)^m \left(1 - \frac{m}{n}\right)^{n-m}$. Using Stirling's formula,

$\binom{n}{m}\left(\frac{m}{n}\right)^m\left(1-\frac{m}{n}\right)^{n-m} \leq \sqrt{\frac{n}{2\pi m(n-m)}} \exp(1/(12n)) \leq \sqrt{\frac{n}{6m(n-m)}}$ for $n \geq 4$. Note that $\frac{n}{6m(n-m)} < \frac{1}{64}$ when $m = 12$ and $n \geq 100$, and then for fixed $n$ the denominator increases as $m$ increases, up to $n/2$, so taking the square root, we see the result holds for the range $12 \leq m \leq (n/2)+1, n \geq 100$.

**Lemma** 4. For $n \geq 100$, $f(m,n)$ is an increasing function of $m$ for $m \leq 11$.

Proof: In fact, we shall show this by estimates only for $n \geq 3200$, and then just numerically check the cases in between. This may seem like laziness, but even the estimate we use is a lot of work, and we were not prepared to go further with the analysis and algebra seemingly needed to do better. We shall use a very direct approach with elementary expansions, and there is room in our method for sharpening the estimates. Another promising approach that perhaps could yield a sharper estimate would be to represent the cumulative binomial with the incomplete beta function, but we won't use that here.

$$f(n,m+1) = \sum_{k \leq m} \binom{n}{k} \left(\frac{m+1}{n}\right)^k \left(1 - \frac{m+1}{n}\right)^{n-k}$$

$$= \sum_{k \leq m-1} \binom{n}{k} \left(\frac{m+1}{n}\right)^k \left(1 - \frac{m+1}{n}\right)^{n-k} + \binom{n}{m}\left(\frac{m+1}{n}\right)^m \left(1 - \frac{m+1}{n}\right)^{n-m}$$

$$= f(n,m) + \sum_{k \leq m-1} \binom{n}{k}\left(\frac{m+1}{n}\right)^k \left(1-\frac{m+1}{n}\right)^{n-k} - \sum_{k \leq m-1} \binom{n}{k}\left(\frac{m}{n}\right)^k \left(1-\frac{m}{n}\right)^{n-k}$$

$$+ \binom{n}{m}\left(\frac{m+1}{n}\right)^m\left(1-\frac{m+1}{n}\right)^{n-m}.$$

Now

$$\left(1-\frac{m}{n}\right)^{n-k} = \exp\left((n-k)\log\left(1-\frac{m}{n}\right)\right) = \exp\left(-m + k\frac{m}{n} - \frac{m^2(n-k)}{2n^2} - \frac{m^3(n-k)}{3n^3} - \ldots\right).$$

Let $r = \frac{m}{n}, r' = \frac{m+1}{n}$. These will be small since we intend to have $n$ large compared to $m$. Then $\left(1-\frac{m}{n}\right)^{n-k} = \exp\left(-m + kr - \frac{m(n-k)}{n}\left(\frac{r}{2} + \frac{r^2}{3} + \ldots\right)\right)$. Similarly, dropping a term $k/n$, $\left(1-\frac{m+1}{n}\right)^{n-k} \geq \exp\left(-(m+1) + kr - \frac{(m+1)(n-k)}{n}\left(\frac{r'}{2} + \frac{r'^2}{3} + \ldots\right)\right)$.

But $(m+1)\left(\frac{r'}{2} + \frac{r'^2}{3} + \ldots\right) - m\left(\frac{r}{2} + \frac{r^2}{3} + \ldots\right) = m\left(\frac{r'-r}{2} + \frac{r'^2-r^2}{3} + \ldots\right) + \frac{r'}{2} + \frac{r'^2}{3} + \ldots$ and

$r'^j - r^j \leq (r'-r) j r'^{(j-1)} = j r'^{(j-1)}/n$ since $r \leq r' = r + 1/n$, so

$$m\left(\frac{r'-r}{2} + \frac{r'^2 - r^2}{3} + ...\right) \leq r\left(\frac{1}{2} + \frac{2r'}{3} + \frac{3r'^2}{4} ...\right) \leq \frac{r}{2}\left(1 + \frac{2r'}{1-r'}\right) = \frac{r(1+r')}{2(1-r')}, \text{ and}$$

$$\frac{r'}{2} + \frac{r'^2}{3} + ... \leq \frac{r'}{2(1-r')}, \text{ so}$$

$$\exp\left(-\frac{(m+1)(n-k)}{n}\left(\frac{r'}{2} + \frac{r'^2}{3} + ...\right) + \frac{m(n-k)}{n}\left(\frac{r}{2} + \frac{r^2}{3} + ...\right)\right) \geq \exp\left(-\frac{r(1+r') + r'}{2(1-r')}\right).$$

Let $\lambda = \exp\left(-\frac{r(1+r') + r'}{2(1-r')}\right)$. So we can now write

$$\left(\frac{m}{n}\right)^k \left(1 - \frac{m}{n}\right)^{n-k} - \left(\frac{m+1}{n}\right)^k \left(1 - \frac{m+1}{n}\right)^{n-k} \leq$$

$$\exp\left(-(m+1) + kr - m(1-r)\left(\frac{r}{2} + \frac{r^2}{3} + ...\right)\right)\left\{e\left(\frac{m}{n}\right)^k - \lambda\left(\frac{m+1}{n}\right)^k\right\}, \text{ and summing,}$$

$$\sum_{k \leq m-1} \binom{n}{k}\left(\frac{m+1}{n}\right)^k \left(1 - \frac{m+1}{n}\right)^{n-k} - \sum_{k \leq m-1} \binom{n}{k}\left(\frac{m}{n}\right)^k \left(1 - \frac{m}{n}\right)^{n-k} \leq$$

$$\exp\left(-(m+1) - m(1-r)\left(\frac{r}{2} + \frac{r^2}{3} + ...\right)\right)\left\{\sum_{k=0}^{m-1} \frac{n!}{n^k} e^{kr}\left(e\frac{m^k}{k!} - \lambda\frac{(m+1)^k}{k!}\right)\right\}. \text{ At this point we}$$

note that it may be shown with standard methods that $\sum_{k=0}^{m-1}\left(e\frac{m^k}{k!} - \frac{(m+1)^k}{k!}\right)$ is less than $\frac{(m+1)^m}{m!}$, so if $n$ were very large so that $r, r'$ and $\delta$ were about zero and $\lambda$ were about one, then we would have $e^{-(m+1)} \sum_{k=0}^{m-1}\left(e\frac{m^k}{k!} - \frac{(m+1)^k}{k!}\right) < e^{-(m+1)} \frac{(m+1)^m}{m!}$, and this right side is about equal to $\binom{n}{m}\left(\frac{m+1}{n}\right)^m \left(1 - \frac{m+1}{n}\right)^{n-m}$, and the result would be proved. So it is just a matter of finding out how big $n$ has to be to make the fudge factors small enough so the inequality holds. It would be possible to use analytical estimates of

$\sum_{k=0}^{m-1}\left(e\frac{m^k}{k!} - \frac{(m+1)^k}{k!}\right)$ to see how much slack we have and estimate the size of $n$ needed

that way, but we are satisfied to just calculate for the finitely many cases we have to find the answer. We shall not try to obtain a really good estimate, because our computer is very fast and a sloppy estimate that leaves lots of cases to be checked is no problem.

Next we estimate in a similar way

$$\binom{n}{m}\left(\frac{m+1}{n}\right)^m\left(1-\frac{m+1}{n}\right)^{n-m} >$$

$$\frac{(m+1)^m}{m!}\frac{n(n-1)\ldots(n-m+1)}{n^m}\exp\left(-(m+1)+(m+1)r-(m+1)(1-r)\left(\frac{r'}{2}+\frac{r'^2}{3}+\ldots\right)\right)$$

Using log expansions of factors $1-1/k$ one may show that
$\frac{n(n-1)\ldots(n-m+1)}{n^m} \geq \exp\left(\frac{-mr(1+r)}{2}\right)$. Then to prove the lemma it suffices to show

$$\exp\left(-(m+1)-m(1-r)\left(\frac{r}{2}+\frac{r^2}{3}+\ldots\right)\right)\left\{\sum_{k=0}^{m-1}e^{kr}\left(e\frac{m^k}{k!}-\lambda\frac{(m+1)^k}{k!}\right)\right\}$$

$$< \frac{(m+1)^m}{m!}\exp\left(-(m+1)-\frac{mr(1+r)}{2}+(m+1)r-(m+1)(1-r)\left(\frac{r'}{2}+\frac{r'^2}{3}+\ldots\right)\right).$$ We have

already noted that $(m+1)\left(\frac{r'}{2}+\frac{r'^2}{3}+\ldots\right)-m\left(\frac{r}{2}+\frac{r^2}{3}+\ldots\right) < \frac{r(1+r')+r'}{2(1-r')}$. So it suffices to show $\sum_{k=0}^{m-1}e^{kr}\left(e\frac{m^k}{k!}-\lambda\frac{(m+1)^k}{k!}\right) \leq \frac{(m+1)^m}{m!}\exp\left(-\frac{mr(1+r)}{2}+(m+1)r-\frac{r(1+r')+r'}{2(1-r')}\right)$. If

$n$ is greater than 120 and $m \leq 11$, then $r' < 1/10$ and $r'-r < 1/120$, and the term being exponentiated on the far right is easily seen to be positive. So our final test is to show

$\sum_{k=0}^{m-1}e^{kr}\left(e\frac{m^k}{k!}-\lambda\frac{(m+1)^k}{k!}\right) \leq \frac{(m+1)^m}{m!}$ for $m \leq 11$ and sufficiently large $n$. Note that if for some $m$, this inequality holds for some $n$, where $r, r'$ and $\lambda$ are determined from $m$ and $n$ as above, then it holds for all larger $n$, because $r$ decreases with $n$ and $\lambda$ increases with $n$. So that's it! By calculation we simply find some value of $n$ so that this inequality holds, for $m \leq 11$.

We report that $n = 3200$ makes the above inequality hold for all $m \leq 11$ (a smaller $n$ works for the smaller $m$, if we cared).

Then for all $m \leq 11$, $100 \leq n < 3200$, we simply verified by calculation that $f(n,m)$ is an increasing function of $m$ (in fact, it is true for smaller $n$ as well, down to around 40).

That concludes our *really ugly* proof of lemma 4, but we can state from experience that getting a sharper estimate in order to save our computer some work is not pretty either.

**Lemma 5**. $\sum_{k \leq m}\binom{n}{k}x^k(1-x)^{n-k}$ is a decreasing function of $x, 0 \leq x \leq 1$.

Proof: This is obvious upon noting that the probability of a binomial random variable being less than some amount decreases as the probability of success in a trial increases. □

So if we were to graph $\sum_{k \leq nx} \binom{n}{k} x^k (1-x)^{n-k}$ as a function of x for $\frac{1}{n} \leq x \leq \frac{1}{2}$, we would see jump discontinuities with positive jumps at the points $x = \frac{m}{n}$, and the heights at the discontinuities increasing with $m$ for a while (depending on $n$), and a decreasing graph between the discontinuities.

**Theorem**. Assume the $X_i$ are identically distributed. $P(\sum_{i=1}^{n} X_i \leq n)$ attains its minimum value $L = \left(1 - \frac{1}{n+1}\right)^n$, with $p_i(0) = 1 - \frac{1}{n+1}$, $p_i(n+1) = \frac{1}{n+1}$, $i = 1,...,n$.

Proof: Let $2 \leq j \leq n, n \geq 100$. Then $\frac{m}{n} \leq \frac{1}{j} < \frac{m+1}{n}$ for some $1 \leq m \leq n/2$, so $m \leq \frac{n}{j} < m+1$. We have $\sum_{kj \leq n} \binom{n}{k} \left(\frac{1}{j}\right)^k \left(1-\frac{1}{j}\right)^{n-k} = \sum_{k \leq m} \binom{n}{k} \left(\frac{1}{j}\right)^k \left(1-\frac{1}{j}\right)^{n-k}$ which by Lemma 5 is greater than or equal $\sum_{k \leq m} \binom{n}{k} \left(\frac{m+1}{n}\right)^k \left(1-\frac{m+1}{n}\right)^{n-k}$. If $m \geq 11$ (note our $m$ here is one greater than in the function of Lemma 3), by Lemma 3 this is greater than 3/8 which is greater than $\left(1-\frac{1}{n+1}\right)^n$. If $1 \leq m \leq 10$, Lemma 4 asserts that the worst case is at $m = 1$, and it is easily checked that $\left(1-\frac{2}{n}\right)^n + n\frac{2}{n}\left(1-\frac{2}{n}\right)^{n-1} = \left(1-\frac{2}{n}\right)^{n-1}\left(3-\frac{2}{n}\right)$ is greater than $\left(1-\frac{1}{n+1}\right)^n$ when $n \geq 100$ (actually much smaller).

That does it, except for the finitely many cases for $2 \leq n \leq 99$, which we checked with our tireless computer. □

## References

[1] Kaas, R. and Buhrman, J. M. (1980). Mean, median and mode in binomial distributions. *Statist. Neerlandica* **34** 13--18.
[2] Winker, P. (2009). Puzzled by Probability and Intuition, *Communications of the ACM*, Vol. 52 no. 8, 104-105.